# The Occurrence-in-subtuple problem

Alexander Platzer

original: November, 2008

corrected version: November, 2014

**Abstract**

As we go along with a bioinformatics analysis, we stumbled over a new combinatorial question. Although the problem is a very special one, there are maybe more applications than only this one we have. This text is mainly about the general combinatorial problem, the exact solution and its derivation. An outline of a real problem of this type is in the discussion.

# 1. The problem

**Definitions**:

**Given** are:
*SN* ... a set of *N* different elements
*SNN* ... a tuple of *N* different elements occurring *N*-1 times
*Sx* ... a subtuple of *SNN*
*Sy* ... a subset of *SN*

**Derived** is:
*Sz* ... the subset of *Sy* of all elements of *Sy* occurring in *Sx*

With the variables:
*N* ... the number of different elements in *SN*
*N-1* ... the count of occurrences of every different element in *SNN* |*SNN*| = *N*(N-1)*
|*Sx*| = *x*
|*Sy*| = *y*
|*Sz*| = *z*

**The question:**
We have *N* different elements occurring *N-1* times. From this tuple we take *x* elements. We have a subset of set *SN* with the number of *y* elements. It is searched in *Sx* elements for the *Sy* elements and there will be *z* elements found. How likely is it to find *z* or more elements of *y* in *x*?
(an illustration of this question is in figure 1)

Or with the general declarations:
What is the ratio of permutations of *SNN,* which results in an equal or higher *z* as given, to all possible permutations of *SNN*?

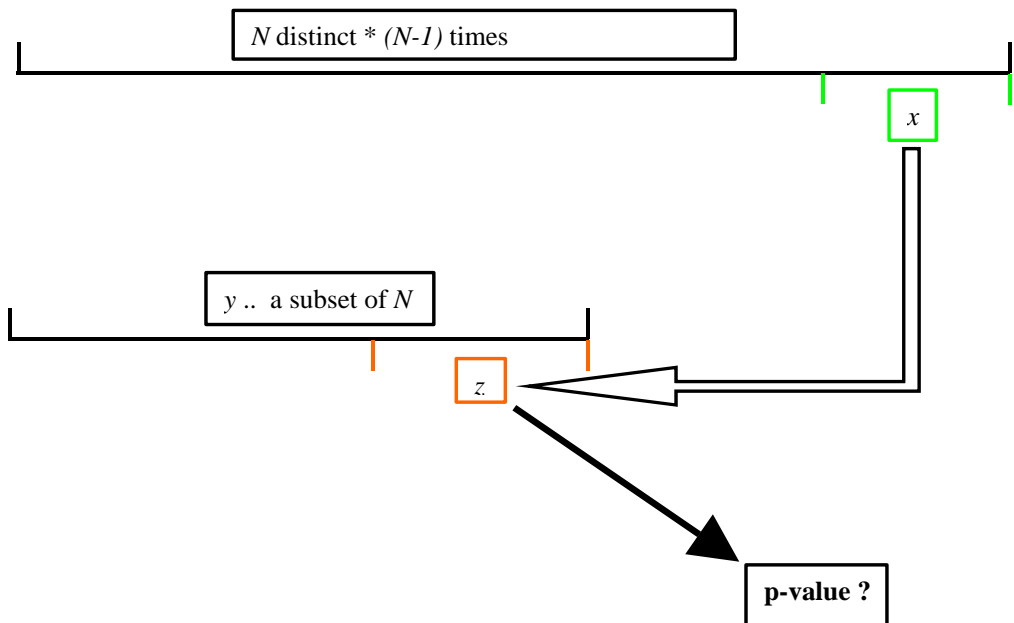

Figure 1. The combinatorial problem illustrated

## 2. Derivation

First of all it is important that the result should be a p-value, so we need the ratio of permutations fulfilling the condition (*z* or more elements of *Sy* are in *Sx*), not the combinations. To get permutations of combinations in *Sx* and the rest the factor *x!\*(N\*(N-1)-x)!* is needed. And at the end we divide by all permutations possible *( (N\*(N-1))!* , the factor becomes

$$\frac{x!*(N*(N-1)-x)!}{(N*(N-1))!} \text{ which is}$$

$$\frac{1}{\binom{N*(N-1)}{x}} \quad (1)$$

if *x=1*, *y=z=1* then this number of combinations is *N-1* (because *N-1* elements of one type exist) if *x=2*, *y=z=1* for 2 times one element of *Sz* in *Sx* the number of combinations is *(N-1)\*(N-2) / 2* ; *(N-1)* possibilities to choose the first and *(N-2)* possibilities to choose the second, but the two chosen elements cannot be distinguished so divide by two) and so on ->

$$\binom{N-1}{how\,many\,times\,one\,element}$$ and for the question, how many combinations exist with one or

more times one element of *Sz* in *Sx* (still *x=2*, *y=z=1*) it is the sum of *N-1* and *(N-1)\*(N-2)* , but notice that this are only the combinations of the elements of Sx which occur in *Sz*, each of this combination has another side (for the elements of *Sx* which are not occurring in *Sz*).
For a number of w elements equal to one element in *Sz* in *Sx*, it can be seen as 2 times k-combinations:

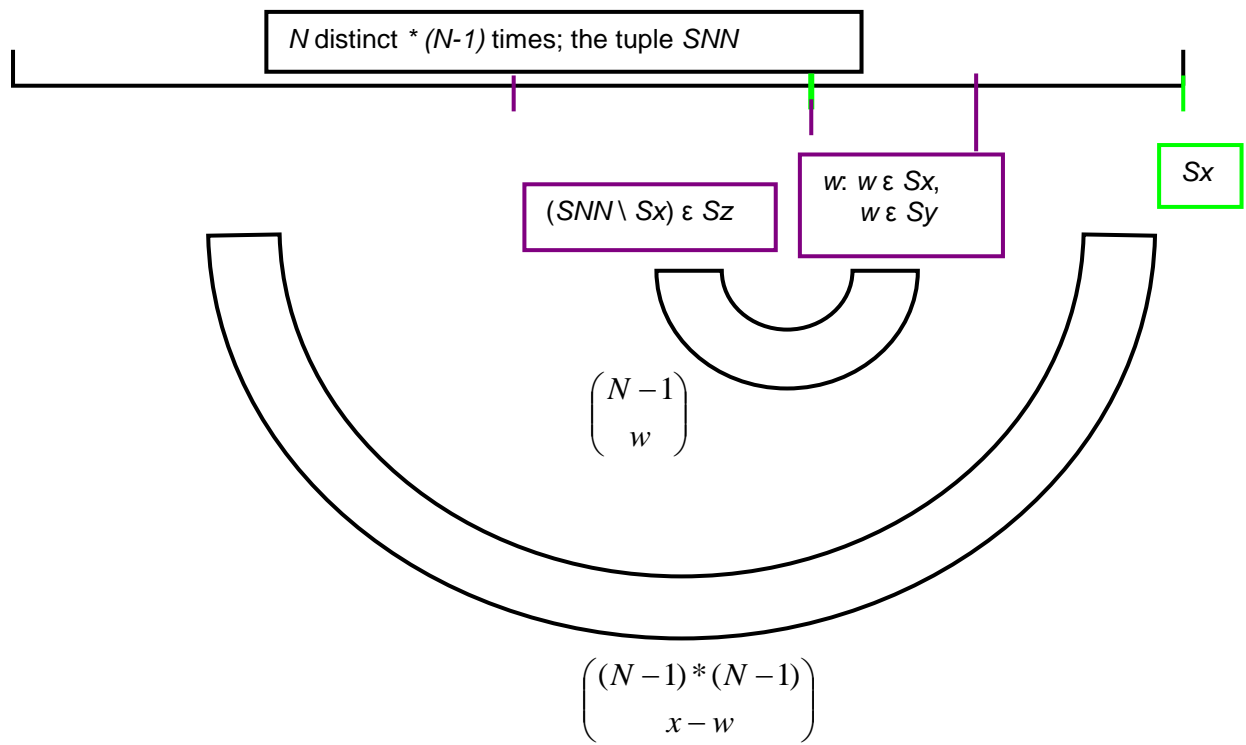

so the count of combinations of w elements in *SNN* element in *Sz* in *Sx*, still *y=z=1*, is

$$\binom{N-1}{w} * \binom{(N-1)*(N-1)}{x-w}$$ , for different *w*'s just summing up.

If we change to more than one *z*, but still *y=z* then *w* is changed to a vector containing the count of every element of *Sz* in *Sx*. Moreover, there must be two constraints, first the sum of w must not be larger than *x* and no component of *w* may be larger than *n-1* (plus as before every component of *w* must be 1 or larger). If we have *k=z* and sum up all combinations (of possible values for each component of *w*):

$$\sum_{w1...k=1}^{w1...k \leq n-1 \wedge \sum_{j=1}^{k} wj \leq x} \left( \left( \prod_{j=1}^{k} \binom{n-1}{w_j} \right) * \binom{(n-y)*(n-1)}{x - \sum_{j=1}^{k} w_j} \right) \quad (2)$$

for the case *y>=z* it is to include that a combination of *z* or more elements are chosen from *Sy*: (vector w is renamed to vector i)

$$\sum_{k=z}^{y} \binom{y}{k} * \sum_{i1...k=1}^{i1...k \leq n-1 \wedge \sum_{j=1}^{k} ij \leq x} \left( \left( \prod_{j=1}^{k} \binom{n-1}{i_j} \right) * \binom{(n-y)*(n-1)}{x - \sum_{j=1}^{k} ij} \right) \quad (3)$$

and with the factor for make permutations out of combinations (1):

$$\frac{\sum_{k=z}^{y} \binom{y}{k} * \sum_{i1...k=1}^{i1...k \leq n-1 \wedge \sum_{j=1}^{k} ij \leq x} \left( \left( \prod_{j=1}^{k} \binom{n-1}{i_j} \right) * \binom{(n-y)*(n-1)}{x - \sum_{j=1}^{k} ij} \right)}{\binom{n*(n-1)}{x}} \quad (4)$$

we have now a formula for the problem, but it can be seen easily that it will need large computing power to solve it with large parameter values. The order of this formula is O(*y \* (n-1)^y*), so we have only changed one combinatorial problem into another (faster, but still incalculable in terms of time). Therefore, the next task is to transform this formula into one with a lower order.
The bad thing in it is the second sum, there the power of y occurs.

First, we expand to:

$$\frac{\sum_{k=z}^{y} \binom{y}{k} * \sum_{s=k}^{x} \left( \sum_{i1...k=1}^{i1...k \leq n-1 \wedge \sum_{j=1}^{k} ij = s} \left( \prod_{j=1}^{k} \binom{n-1}{i_j} \right) * \binom{(n-y)*(n-1)}{x-s} \right)}{\binom{n*(n-1)}{x}}$$

(5)

One thing to remark is that the restriction $i_{1...k} \leq n-1$ is only there for information, it is fulfilled of the formula alone because the binomial coefficient of $\binom{N-1}{a}$ with *a > N-1* is 0.

Now we look closer to

$$\sum_{i1...k=1}^{i1...k \leq n-1 \wedge \sum_{j=1}^{k} ij = s} \left( \prod_{j=1}^{k} \binom{n-1}{i_j} \right) \quad (6)$$

If the vector components would not start with 1 but with 0 it would look a trifle better:

$$\sum_{i_{1\ldots k}=0}^{i_{1\ldots k}\leq n-1 \wedge \sum_{j=1}^{k}i_j=s} (\prod_{j=1}^{k}\binom{n-1}{i_j}) \qquad (7)$$

(Caution: In the next steps in between, the variable names are not always the same as before)

This we can transform with a generalization of the Vandermonde's identity.

Vandermonde's identity: $\sum_{j}\binom{m}{j}*\binom{n-m}{k-j}=\binom{n}{k}$ \qquad (8)

Generalization: $\sum_{k_{1\ldots y}=0}^{k_{1\ldots y}<=n}\binom{n}{k_1}*\binom{n}{k_2}*\binom{n}{k_3}*\ldots*\binom{n}{x-\sum_{j=1}^{y}k_j}=\binom{(y+1)*n}{x}$ \qquad (9)

(generalization of Vandermonde's identity is similar to its algebraic proof (2008), only with *k+1* polynomials instead of 2)

Here is to see that the conditions
- The sum of all lower parts of the binomial coefficients must be *x* (because of the last term with $x-\sum_{j=1}^{y}k_j$ )
- if a lower part of the binomial coefficients is larger than *n* or lower than zero then the factor is zero and the corresponding addend is also zero (this holds also for $x-\sum_{j=1}^{y}k_j$ )

With that we can transform (7) into $\binom{k*(n-1)}{s}$. This would be much faster for computation, but we have had (6), not (7).

$$\sum_{i_{1\ldots k}=0}^{i_{1\ldots k}\leq n-1 \wedge \sum_{j=1}^{k}i_j=s} ((\prod_{j=1}^{k}\binom{n-1}{i_j}) \ (7) \ > \ \sum_{i_{1\ldots k}=1}^{i_{1\ldots k}\leq n-1 \wedge \sum_{j=1}^{k}i_j=s} ((\prod_{j=1}^{k}\binom{n-1}{i_j}) \ (6)$$ so we can use the first

formula and subtract the difference to the second.
The difference is somehow a similar problem because the lower part of the last binomial-coefficient $i_k$ is zero if the sum of $i_{1\ldots k-1}$ is *s* (important in (9)). This is the same problem as before but this time with *k-1*.

However, the difference (7) - (6) is not: $k * \sum_{i_{1\ldots k-1}=0}^{i_{1\ldots k-1}\leq n-1 \wedge \sum_{j=1}^{k-1}i_j=s} ((\prod_{j=1}^{k}\binom{n-1}{i_j})$ \qquad (10) because some combinations will then be counted more than once.

To illustrate this:

We have the sum for every combination of *i* (the components of vector *i* are every combination of one value per column):

| 0 | 0 | 0 | 0 | 0 |
|---|---|---|---|---|
| 1 | 1 | 1 | 1 | 1 |
| ... | ... | ... | ... | ... |
| n-1 | n-1 | n-1 | n-1 | n-1 |

but we need only:

| 1 | 1 | 1 | 1 | 1 |
|---|---|---|---|---|
| ... | ... | ... | ... | ... |
| n-1 | n-1 | n-1 | n-1 | n-1 |

(Note that the constraint of the fixed sum of all components of *i* is fulfilled with a term $s - \sum_{j=1}^{k-1} k_j$, so the tables above have *k-1* colums)

Therefore, we have to substract every combination of *i* with at least one component is zero.

However, the result of $k * \sum_{\substack{i1...k-1 = 0}}^{i1...k-1 \leq n-1 \wedge \sum_{j=1}^{k-1} i_j = s} ((\prod_{j=1}^{k} \binom{n-1}{i_j}))$ *(10)* would be:

| 0 | 0 | 0 | 0 | 0 | | 0 | 0 | 0 | 0 | 0 | | 0 | 0 | 0 | 0 | 0 |
|---|---|---|---|---|---|---|---|---|---|---|---|---|---|---|---|---|
| 1 | 1 | 1 | 1 |   | | 1 |   | 1 | 1 | 1 | | 1 | 1 |   | 1 | 1 |
| ... | ... | ... | ... |   | | ... |   | ... | ... | ... | | ... | ... |   | ... | ... |
| n-1 | n-1 | n-1 | n-1 |   | | n-1 |   | n-1 | n-1 | n-1 | | n-1 | n-1 |   | n-1 | n-1 |

.....

Then we count every combination of two times 0 more than once, so we have add again these combinations:

| 0 | 0 | 0 | 0 | 0 | | 0 | 0 | 0 | 0 | 0 | | 0 | 0 | 0 | 0 | 0 |
|---|---|---|---|---|---|---|---|---|---|---|---|---|---|---|---|---|
|   | 1 | 1 | 1 |   | |   | 1 |   | 1 | 1 | |   | 1 | 1 |   | 1 |
|   | ... | ... | ... |   | |   | ... |   | ... | ... | |   | ... | ... |   | ... |
|   | n-1 | n-1 | n-1 |   | |   | n-1 |   | n-1 | n-1 | |   | n-1 | n-1 |   | n-1 |

..... Here we have once again the similar problem with *k-2*, now $\binom{k}{2}$ times. If we go further and further with alternating add and subtract we end at $\binom{k}{k}$ where all combinations are counted properly ->

$$\sum_{\substack{i1...k = 1}}^{i1...k \leq n-1 \wedge \sum_{j=1}^{k} i_j = s} ((\prod_{j=1}^{k} \binom{n-1}{i_j})) \qquad (6)$$

=

$$\binom{k*(n-1)}{s} + \sum_{j=1}^{k-1} (-1)^j * \binom{k}{j} * \binom{(k-j)*(n-1)}{s} \qquad (11)$$

(11) combined with (5) results in:

$$\frac{\sum_{k=z}^{y}\binom{y}{k}*\sum_{s=k}^{x}(\binom{k*(n-1)}{s}+\sum_{j=1}^{k-1}(-1)^{j}*\binom{k}{j}*\binom{(k-j)*(n-1)}{s})*\binom{(n-y)*(n-1)}{x-s}}{\binom{n*(n-1)}{x}}$$

with the order O($y * x^2$) and because the numbers are growing larger (but still needed to be computed as full integers) with a little term to the power of 4.

Because of the long drawn out derivation, the formula (12) was tested against (4) in a not too time-consuming parameter range (*N* < 20 and *y* < 20). The formula (4) again was tested against the full amount of permutations, but only *N* < 5 was feasible.

For the author of this derivation it seems remarkable that the constraints (which are time-consuming) e.g. $i(1)...(k) \leq n-1 \wedge \sum_{j=1}^{k} i(j) = s$ can be led back to the constraint of factorial(lower than 0) = 0.

## 3. Discussion

Even though the order of O($x^3$) does not look so bad, there should be a faster approximate formula with a guaranteed bound of error. However, here the motivation was 'let us try to see if it is possible to solve it exactly'.
The first occurrence/use of this type of problem was in gene expression data (as far as we know).
Assume we have n elements which are somehow (but unknown) regulating each other. So there are n*n possible regulations or n*(n-1) regulations without elements regulating itself. Of this n*(n-1) possible regulations we could take x regulations for some reason (e.g. it is the result of a method to filter all possible regulations). Additionally given is a list of y elements, which are known as regulators. (if we annotate the regulations as lines of 'A -> B' then the regulators are just the elements on the left side).
So we can search in (the first column of) the x chosen regulations for this y elements and we will find a count of z elements.
Now the question arises, how likely is this? Alternatively, the other way round, is it just by chance to get z of y elements in our filtered regulations?

Such a scenario is not quite frequent and it sounds strange to know that something is a regulator without knowing what is exactly regulated by it. Nevertheless, it is possible to know that a gene very likely is a transcription factor, even if it is not quite sure which other genes are exactly regulated of it.